\documentclass[12pt,letterpaper,onecolumn]{IEEEtran}
\usepackage{amsmath}
\usepackage{subfigure}
\usepackage{graphicx}
\graphicspath{{EPS/}}

\title{Network Kriging}
\author{David B.\ Chua, Eric D.\ Kolaczyk, Mark Crovella%
\thanks{David B.\ Chua (dchua@math.bu.edu) and Eric D.\ Kolaczyk
  (kolaczyk@math.bu.edu) are with the Dept.\ of Mathematics and Statistics
  at Boston University. Mark Crovella (crovella@cs.bu.edu) is with the
  Computer Science Dept.\ at Boston University. 
Part of this work was performed
        while E. Kolaczyk was with the LIAFA group at l'Universit\'e
        de Paris-7, with support from the CNRS, and while M. Crovella was 
  at the Laboratoire d'Informatique de Paris-6 (LIP6), 
  with support from CNRS and Sprint Labs.
        This work was supported in part by NSF
        grants ANI-9986397 and CCR-0325701, and by ONR award N000140310043.}
}

\newcommand\onefig{3in}
\newcommand\twofig{3in}
\newcommand\G{\mathcal{G}}
\newcommand\V{\mathcal{V}}
\newcommand\E{\mathcal{E}}
\newcommand{\calS}{\mathcal{S}}

\newcommand{\Ss}{\calS_s}
\renewcommand\P{\mathcal{P}} 
\newcommand\abs[1]{\lvert #1 \rvert}
\newcommand\R{\mathbf{R}}
\newcommand\Vss{V_{ss}}
\newcommand\Vrs{V_{rs}}
\newcommand\Vsr{V_{sr}}
\newcommand\Vrr{V_{rr}}
\newcommand\tildeGs{\tilde{G}_s}
\newcommand\Gs{G_s}

\newcommand\Gr{G_r}
\newcommand\Bs{B_s}

\newcommand\ahat{\hat{a}}
\newcommand\muhat{\hat{\mu}}
\newcommand\ie{\emph{i.e.}}
\newcommand\eg{\emph{e.g.}}

\newcommand{\Norm}[1]{\left\lVert #1 \right\rVert}
\newcommand{\norm}[1]{\lVert #1 \rVert}
\DeclareMathOperator\bias{Bias}
\DeclareMathOperator\mspe{MSPE}
\DeclareMathOperator\rank{Rank}
\DeclareMathOperator\nul{Null}

\DeclareMathOperator\row{Row}
\DeclareMathOperator\diam{diam}
\newtheorem{theorem}{Theorem}[section]
\newtheorem{proposition}[theorem]{Proposition}

\begin{document}

\maketitle

\begin{abstract}
  Network service providers and customers are often concerned with
  aggregate performance measures that span multiple network paths.
  Unfortunately, forming such network-wide measures can be difficult,
  due to the issues of scale involved. In particular, the number of
  paths grows too rapidly with the number of endpoints to make
  exhaustive measurement practical.  As a result, it is of interest to
  explore the feasibility of methods that dramatically reduce the
  number of paths measured in such situations while maintaining
  acceptable accuracy.

  We cast the problem as one of statistical prediction---in the
  spirit of the so-called `kriging' problem in spatial
  statistics---and show that end-to-end network properties may be
  accurately predicted in many cases using a surprisingly small
  set of carefully chosen paths.  More precisely, we formulate a
  general framework for the prediction problem, propose a class
  of linear predictors for standard quantities of interest (\eg,
  averages, totals, differences) and show that linear algebraic
  methods of subset selection may be used to effectively choose
  which paths to measure. We characterize the performance of the
  resulting methods, both analytically and numerically.  The
  success of our methods derives from the low effective rank of
  routing matrices as encountered in practice, which appears to be
  a new observation in its own right with potentially broad
  implications on network measurement generally.
\end{abstract}

\section{Introduction}
\label{sec:introduction}

In many situations it is important to obtain a network-wide view
of path metrics, such as latency and packet loss rate.  For
example, in overlay networks regular measurement of path
properties is used to select alternate routes.  At the IP level,
path property measurements can be used to monitor network health,
assess user experience, and choose between alternate providers,
among other applications.  Typical examples of systems performing
such measurements include the NLANR AMP project, the RIPE
Test-Traffic Project, and the Internet End-to-end Performance
Monitoring project \cite{NLANR,RIPETT,iepm}.

Unfortunately extending such efforts to large networks can be
difficult, because the number of network paths grows as the square
of the number of network endpoints.  Initial work in this area has
found that it is possible to reduce the number of end-to-end
measurements to the number of ``virtual links'' (identifiable link
subsets)---which typically grows more slowly than the number
of paths---and yet still recover the complete set of end-to-end
path properties exactly \cite{chen03:overlay,shavitt01:computing}.

This quantity is a sharp lower limit that stems from a linear
algebraic analysis of the rank of routing matrices. Measuring even
one path fewer requires one to consider approximations instead of
exact reconstructions.  Specifically, one is faced with the task of 
measuring some paths and predicting the characteristics of others.
The prediction of population characteristics from those of a sample
is a classical problem in the statistical literature.  The most 
well-known version of the predication problem
is perhaps that which occurs in the spatial
sciences, under the name of \emph{kriging} \cite{cressie93:spatial},
where, for example, measurements are taken at series of spatially
distributed wells to enable prediction of oil concentrations
throughout the underlying substrate.

In this paper, we develop a framework for what we term
\emph{network kriging}, the prediction of network path characteristics
based on a small sample. Our methods exploit an observed tendency
in real networks for sharing certain edges between many paths, i.e., a
sort of ``squeezing'' of paths over these edges.
We begin with a discussion of this sharing and the reduced
effective rank of routing matrices in
Section~\ref{sec:routing-matrices}, followed by a development of
our statistical framework and a path selection algorithm in
Section~\ref{sec:prediction-E2E-properties}. In
Section~\ref{sec:empirical-evaluation} we examine the performance
of our methods using delay data from the backbone of the Abilene
network. In Section~\ref{sec:discussion} we conclude with a brief
discussion.

\section{Routing Matrices: Rank versus Effective Rank}
\label{sec:routing-matrices}

\subsection{Background}
\label{sec:background}

We begin by establishing some relevant notation and definitions.
Let $\G = (\V, \E)$ be a strongly connected directed graph, where
the nodes in $\mathcal{V}$ represent network devices and the edges
in $\mathcal{E}$ represent links between those devices.
Additionally, let $\mathcal{P}$ be the set of all paths in the
network, and let $n_v = \abs{\V}$, $n_e = \abs{\E}$, $n_p =
\abs{\P}$ denote respectively the number of devices, links and
paths.
Finally, let $y\in\R^{n_p}$ be the values of a metric measurable 
on all paths $i\in\P$, which is assumed to be a linear function
of the values of the same metric on the edges $j\in\E$, expressed
as $x\in\R^{n_e}$.  We are interested in particular
in the case where $n_p\gg n_e$ and
the linear relation between $y$ and $x$ is given by $y=Gx$, where
$G\in\{0,1\}^{n_p\times n_e}$ is a routing matrix whose entries
simply indicate the traversal of a given link by a given path via
\begin{equation}
  \label{eq:1}
  G_{i,j} = 
  \begin{cases}
    1 & \text{if path $i$ traverses link $j$,} \\
    0 & \text{otherwise.}
  \end{cases}
\end{equation}
For example, if we let $x$ denote the delay times for edges in the
network and let $y$ denote the delay times for paths in the
network, then $y=Gx$.  Additionally, the same relation holds for
$\log(1-\text{loss rate})$.

As explained in Section~\ref{sec:introduction}, our interest in
this paper focuses on the problem of monitoring global
network properties via measurements on some small subset of the
paths.  Note that the question of which paths to monitor is
equivalent to the selection of an appropriate subset of 
rows in $G$, due to the relation $y=Gx$.  Exploiting this insight,
earlier work by Chen and colleagues \cite{chen03:overlay}
shows that in fact one can measure as few as $k^* \sim O(n_e)$
paths and still recover {\it exact} knowledge of all network path
behaviors.\footnote{In \cite{chen04:algebraic} Chen and colleagues
  show that the number $k^*$ of paths needed for their method
  scales at worst like $O(n_v\log n_v)$ in a collection of real
  and simulated networks and they argue that this behavior is to
  be expected in internet networks, due to the high degree of
  sharing between paths that traverse the dense core.}  Their
argument is essentially linear algebraic in nature, and is based
upon the fact that a subset, say $\tilde G$, of only
$k^*=\rank(G)$ independent rows of $G$ are sufficient to span the
range of $G$, \ie, to span the set $\{y\in\R^{n_p} : y=Gx,
x\in\R^{n_e}\}$.  As a result, given the measurements for paths
corresponding to the rows of such a $\tilde G$, measurements for
all other paths may be obtained as a function thereof.  Similar
work may be found in~\cite{nguyen04:link_failures}, in the context
of Boolean algebras, for the problem of detecting link failures.

\subsection{Reduced Rank}
\label{sec:reduced-rank}

Critical to the success of the methodology we propose
in Section~\ref{sec:prediction-E2E-properties} is the concept of 
\emph{effective rank}---the number of independent rows required
to approximate a given matrix to a pre-determined level of tolerance.
Effective rank is an important tool in numerical analysis and
scientific computing (\eg, \cite{golub89}), where it is often used
to reduce the dimensionality of a linear system, generally with the
goal of improving numerical stability.  Here we use it to effect
a significant reduction in measurement requirements above and beyond
the levels achievable through the methodology 
in~\cite{chen03:overlay,chen04:algebraic}.  In particular, we have
found that routing matrices $G$ have an effective rank much smaller
than their actual rank, and as a result a surprisingly small
number of rows are sufficient to adequately approximate the span of $G$.

\begin{figure}[tbp]
  \centering
  \subfigure[Map of Abilene]{%
    \label{fig:abilene-map}
    \includegraphics[width=\twofig]{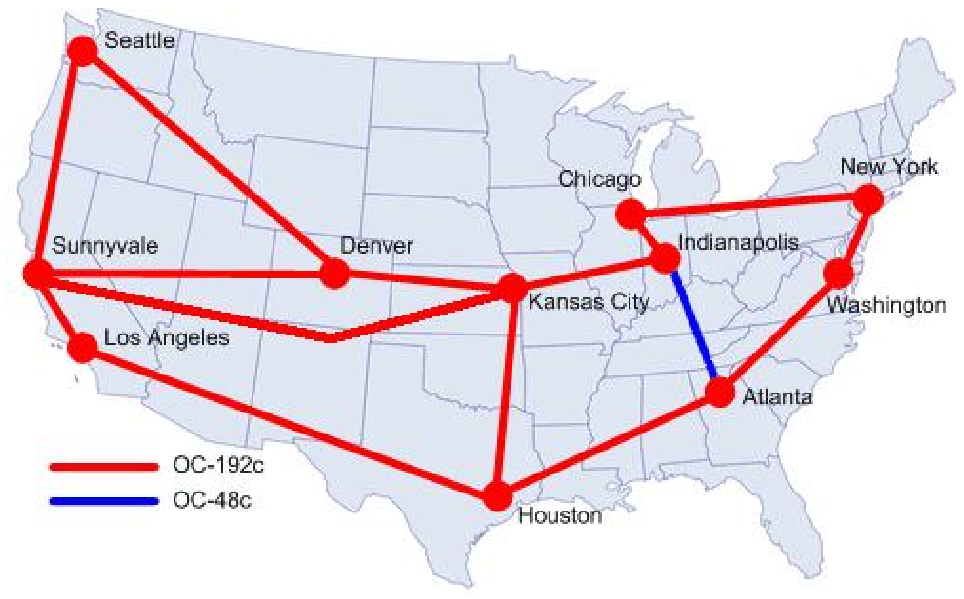}}
  \quad
  \subfigure[Eigenspectrum of $G$]{%
    \label{fig:abilene-spec}
    \includegraphics[width=\twofig]{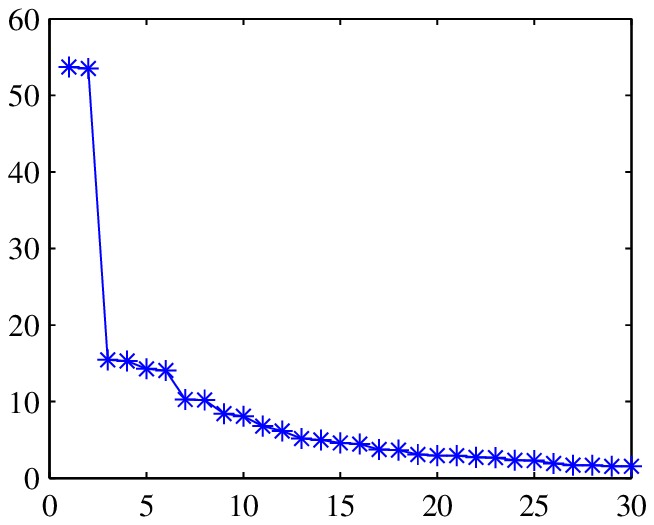}
  }
  \caption{Map of the Abilene network and the eigenspectrum of one
    of its routing matrices.}
  \label{fig:map-and-spec}
\end{figure}

As an illustration, consider the Abilene network shown in
Figure~\ref{fig:abilene-map}; this is a high-performance network
that serves Internet2 (the U.S.\ national research and education
backbone).  The network can be seen to consist of $11$ nodes, but
only $2\times 15=30$ directed links.  Accordingly, a large amount
of sharing of these links can be expected between the $11\times
10=110$ paths on the network. Such sharing would mean a great deal
of similarity between paths and thus fewer ``unique'' paths to
measure. Furthermore, similarities between paths would mean
similarities between the rows of $G$ which suggests that $G$ may
have an effective rank less than 30.

A standard tool for assessing dimensionality and effective rank is
the singular value decomposition (SVD) of $G$, which can be derived from
an eigen-analysis of the matrix $B=G^TG$.  The
eigenvalues of $B$ (\ie, the squares of the singular values of
$G$) are plotted in Figure~\ref{fig:abilene-spec}. The large gap
between the second and third eigenvalues, and the resulting knee
in the spectrum, is evidence of a non-trivial amount of linear
dependence among the rows of $G$. The effective rank of a matrix
is determined by looking for a large gap in the spectrum that
partitions the spectrum into large and small values. Thus the gap
in Figure~\ref{fig:abilene-spec} suggests that as few as two paths
may be sufficient to recover useful information about $y$
in the Abilene network.

\begin{figure}[tbp]
  \centering
  \includegraphics[width=\onefig]{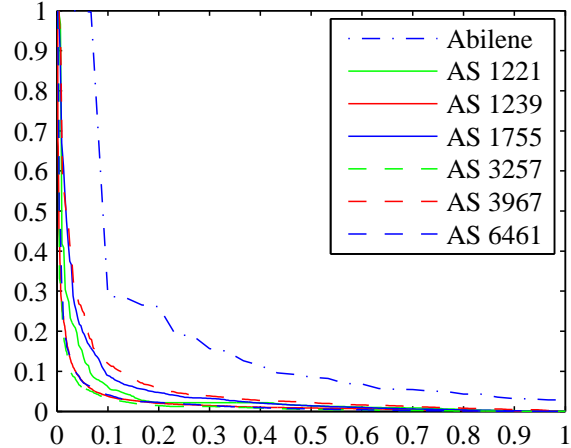}
  \caption{Spectra of $G$ for six networks mapped by the Rocketfuel project and
    Abilene. Note that the spectrum for each network has been re-scaled
    by the largest eigenvalue and the indices have been re-scaled
    to the unit interval. So on the horizontal axis, 1
    corresponds to the $\rank(G)$-th eigenvalue.}
  \label{fig:spec-G}
\end{figure}

Such strong spectral decay appears to be a common property of
routing matrices.  In Figure~\ref{fig:spec-G} we plot the spectra
for five of the networks mapped by the Rocketfuel project~\cite{rocketfuel}. 
The sharp knee that occurs roughly
20\% of the way through is evidence that the effective rank of
these routing matrices is much smaller then their actual rank.
Furthermore, a remarkable amount of similarity can be seen in the
decay of the five spectra.

\begin{figure*}[tbp]
  \vspace*{-0.75in}
  \centerline{
    \includegraphics[width=0.25\textwidth]{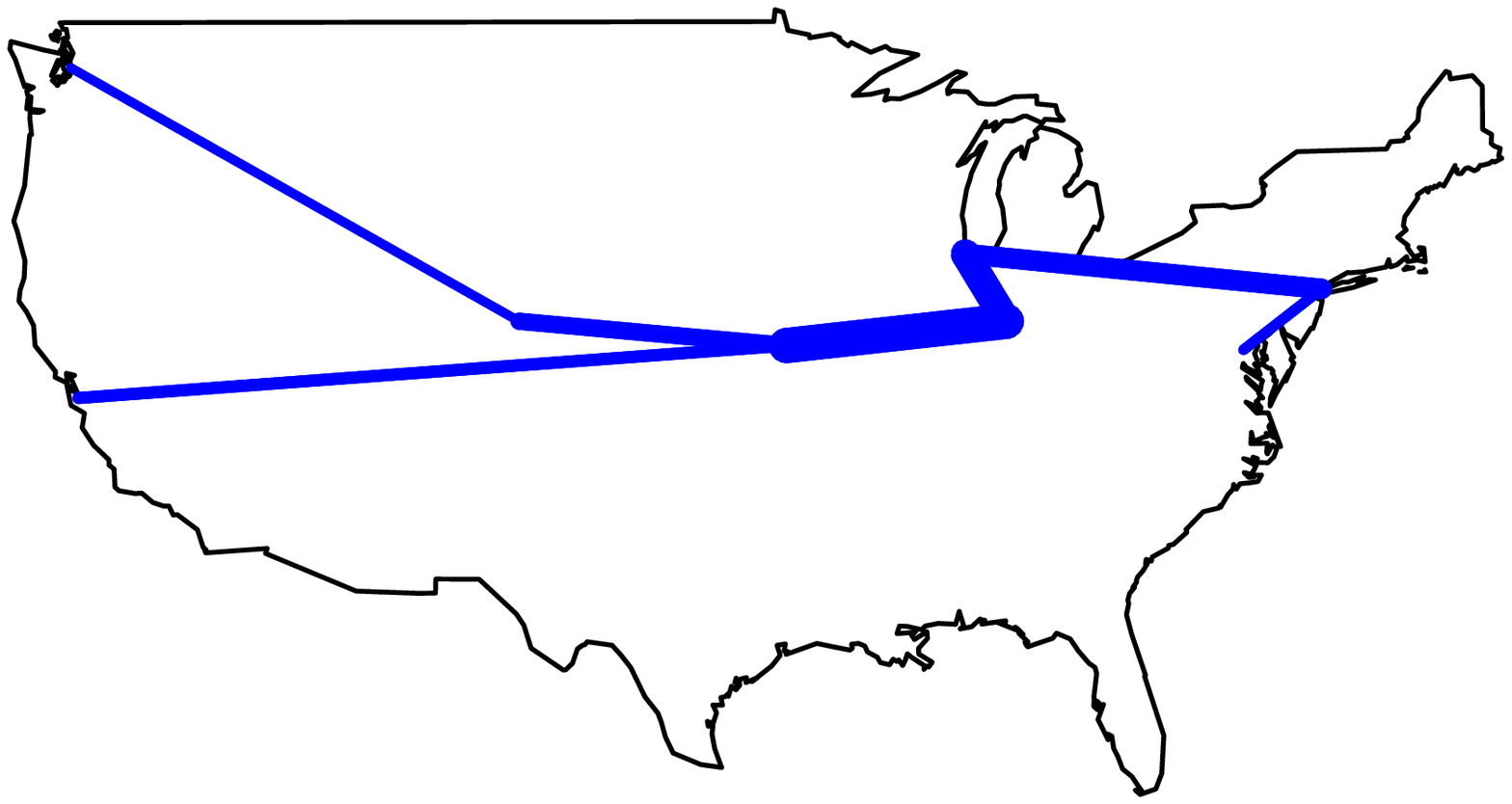}
    \includegraphics[width=0.25\textwidth]{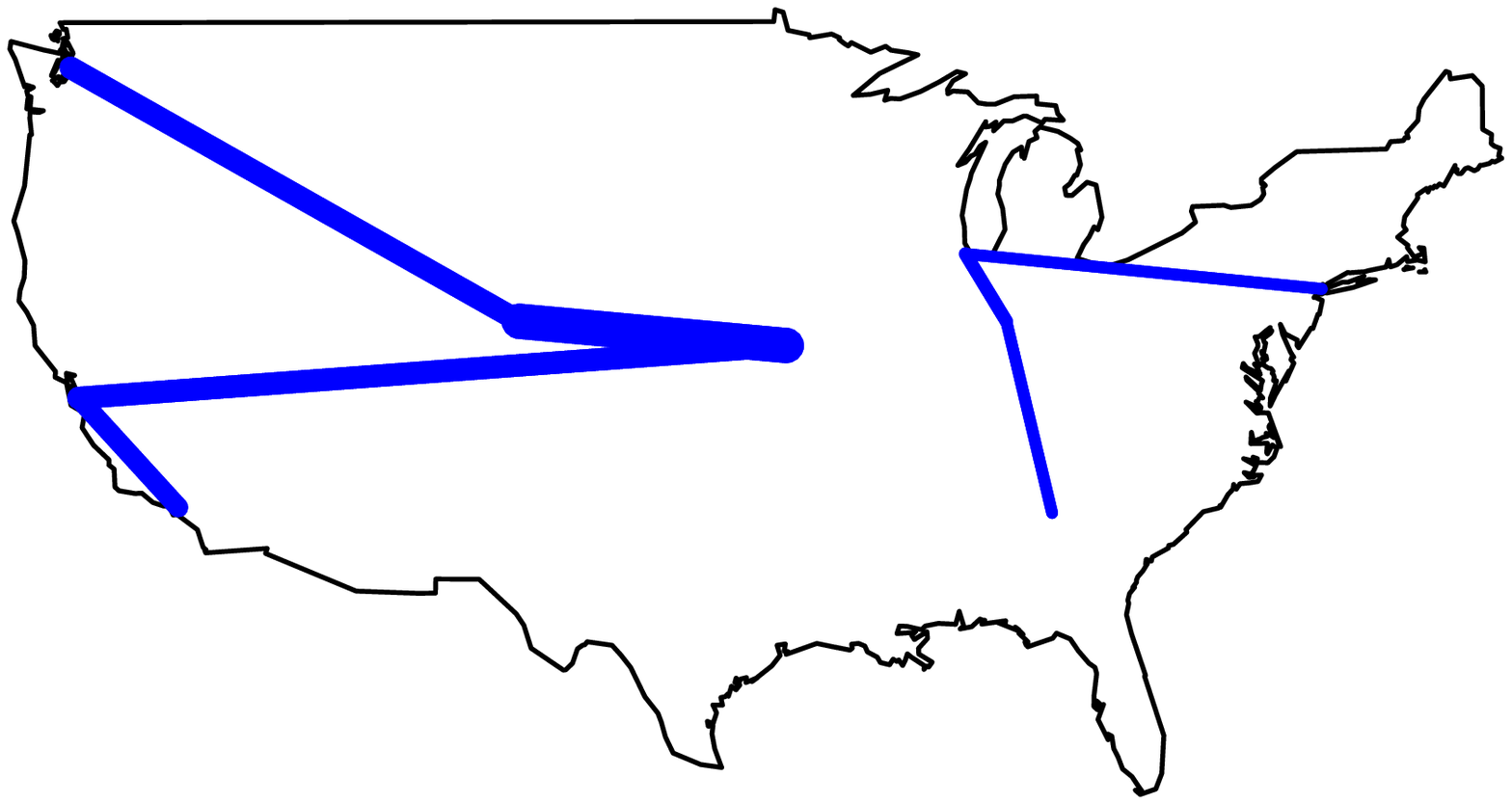}
    \includegraphics[width=0.25\textwidth]{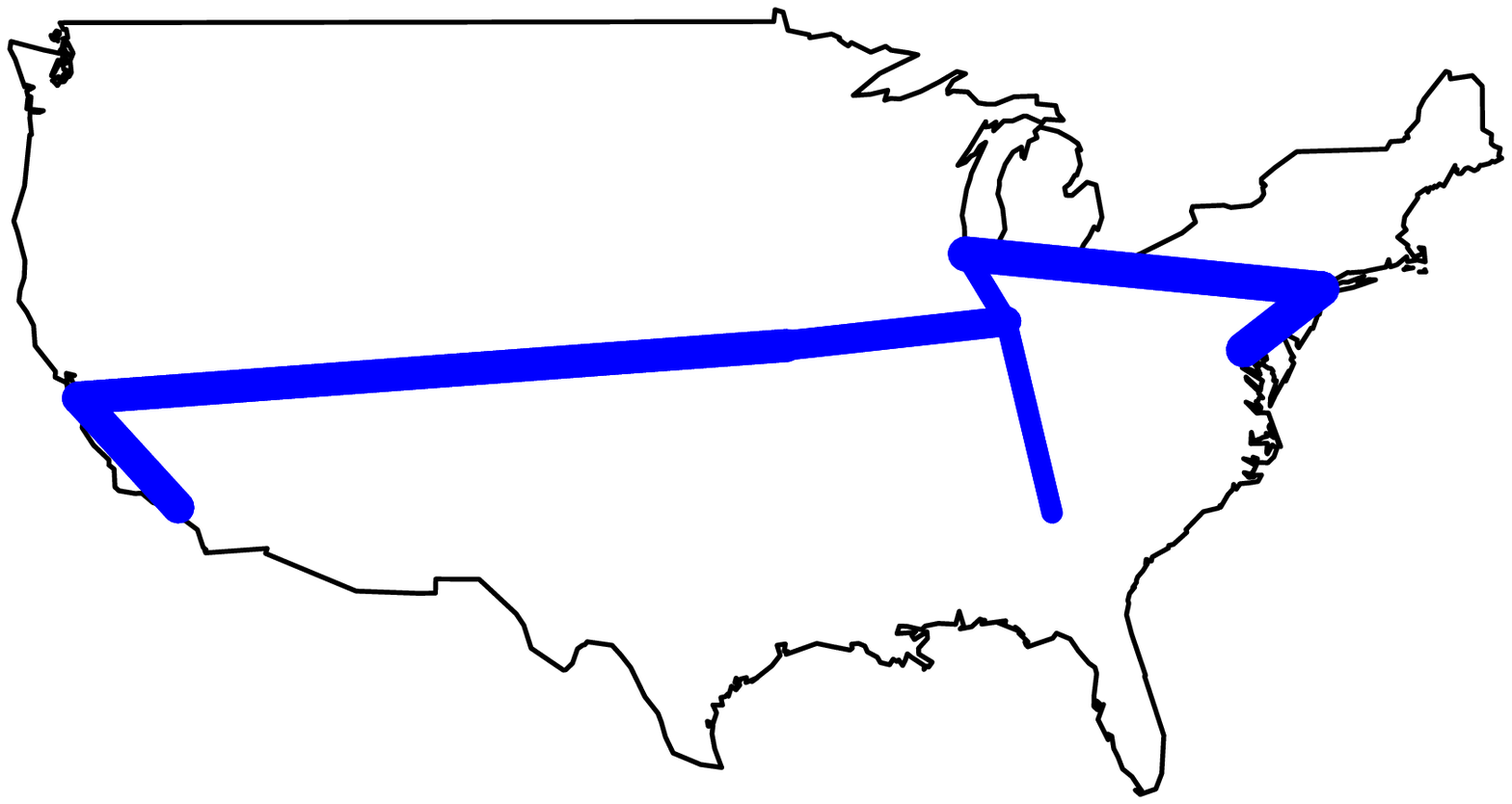}
    \includegraphics[width=0.25\textwidth]{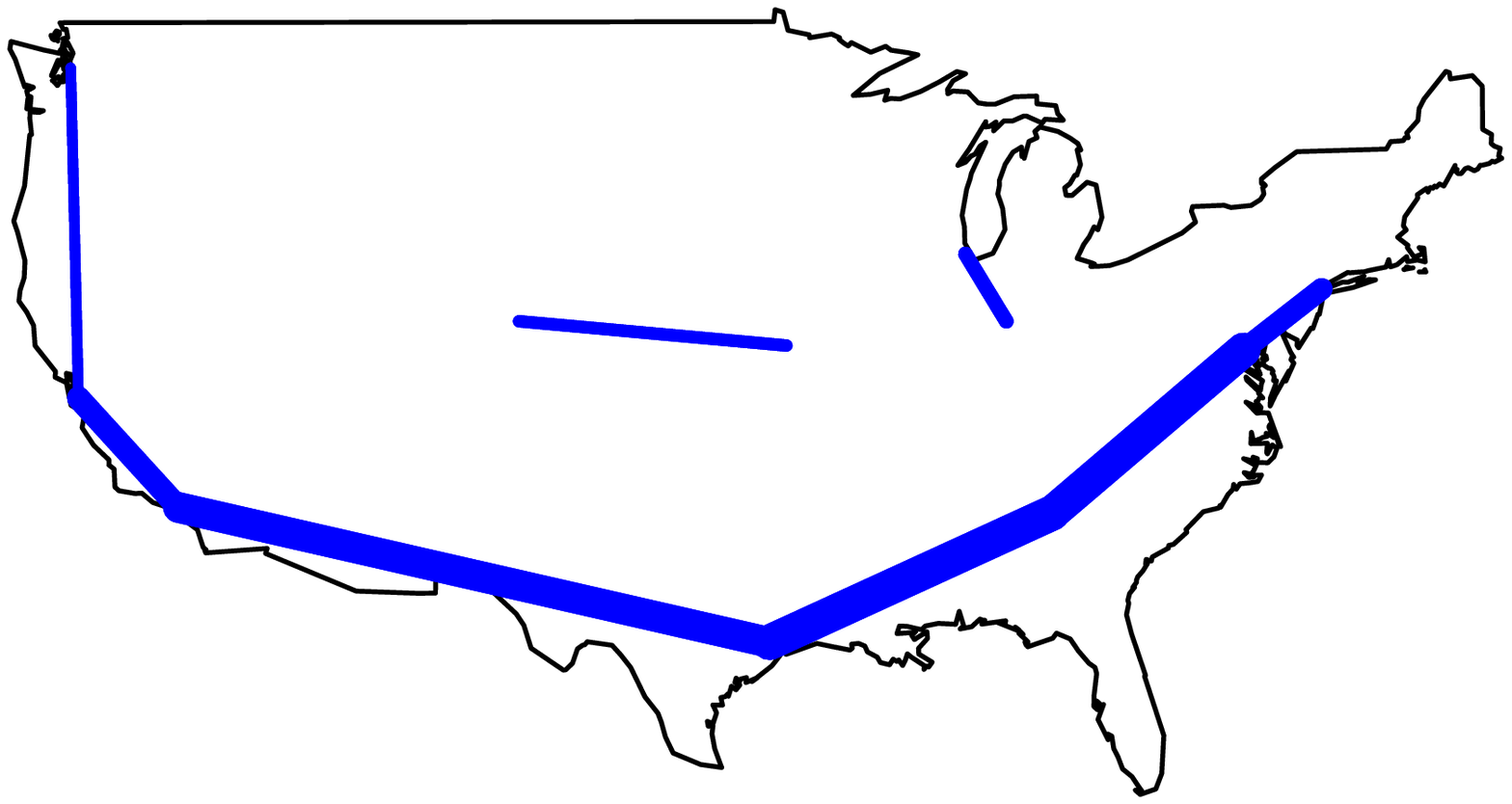}
  }
  \vspace*{-0.75in}
  \centerline{\hfill(a)\hfill\hfill(b)\hfill\hfill(c)\hfill\hfill(d)\hfill}
  \caption{First four distinct eigenvectors of $B=G^TG$ for
    Abilene. Each link is drawn with a thickness that is roughly
    proportional to the magnitude of its corresponding eigenvector
    component. 
  }
  \label{fig:abilene-eigenroutes}
\end{figure*}

To better appreciate the connection between this spectral behavior
and network path properties, we turn to the eigenvectors, which here
may be viewed as orthogonal vectors that capture independent
``patterns'' that occur among the paths in $G$.  For example, the first
eigenvector corresponds to the ``direction'' in link space that
maximizes the path volume of $y$, \ie, $v_1 = \arg
\max_{\norm{x}=1} x^T G^T G x = \arg \max_{\norm{x} = 1} y^Ty$. As
can be seen in Figure~\ref{fig:abilene-eigenroutes}, the energy of
the first eigenvector for Abilene is concentrated along an east-west 
``path'' across the northern part of the network,
with the greatest concentration of energy
occurring at the centrally located Indianapolis-Kansas City link.
The subsequent eigenvectors can be seen to bring in successive
refinements to this picture, with the second and third eigenvectors
further emphasizing connections to and within the ``path'' indicated
by the first, while with the fourth eigenvector we begin to see
evidence of a second east-west ``path'' along the southern part of
the network.  See \cite{chua05:efficient} for additional discussion.

\subsection{Connection to Betweenness}
\label{sec:connection-between}

The evidence in Figures~\ref{fig:abilene-spec} and~\ref{fig:spec-G}
indicates that the effective rank of routing matrices in real 
networks can be noticeably lower than the actual rank.  In the
next section we will show how this phenomenon allows for 
substantial savings in measurement load for the particular problem of
path monitoring that we have chosen to study.  But we believe that the 
implications are in fact much broader.  This would suggest that the issue 
of reduced rank is a topic worth better understanding in and of itself.
For example, from a practical perspective, it would be useful 
to understand how decisions in network design and route management
affect the relative change from actual to effective rank.  In this
regard, connections between the spectral decay of $G$, on the one hand,
and known metrics of topological structure, on the other hand, are
likely to be useful.  While a comprehensive study of this sort is beyond
the scope of the present paper, we describe here a result establishing 
a connection with one such metric.

Specifically, note that the plots in Figure~\ref{fig:abilene-eigenroutes}
confirm our original intuitive notion that the effective rank of $G$
bears an intimate connection with the disproportionate role
played by some links over others in the routing of paths within the
network.  This observation suggests the relevance of the concept of
{\it betweenness centrality}, a concept fundamental in the literature
on social networks \cite{wasserman94} and more recently 
being used in the study of complex networks in the statistical physics 
literature (\eg, \cite{barthelemy04:between,newman04:community}).  
Essentially, betweenness 
centrality measures the number of paths that utilize a specified node,
in the case of `vertex centrality', or a specified link, in the
case of `edge centrality'.\footnote{The raw number of paths utilizing 
a node/edge is typically used under the assumption of unique shortest
path routing; when multiple shortest paths exist, various methods of 
weighted counting have been proposed.  See~\cite{wasserman94}.}

Note that the diagonal elements of the matrix $B=G^TG$ are precisely 
the number of paths routed over their respective links and hence a
measure of the centrality of each link in the network.  The off-diagonal
elements measure the number of paths routed simultaneously over pairs
of links, and might therefore be termed a measure of edge
`co-centrality' or `co-betweenness'.  The co-betweenness $B_{i,j}$
of any two edges $i$ and $j$ will always be bounded above by the 
smaller of the two edge betweenness', 
\ie, $B_{i,j}\le \min\{B_{i,i}, B_{j,j}\}$.  Hence, it might not
be unreasonable to expect that the behavior of the 
eigen-spectrum of $B$ may be related to that of its diagonal, as
the following result shows.
\begin{proposition}
  \label{prop:1}
  Let $B=G^TG$ and, without loss of generality, assume that the
  edges have been ordered so that 
  $B_{1,1}\geq \dots \geq B_{n_e,n_e}$.
  Then, for $k=1,\dots,n$, we have
  $\lambda_k \leq B_{k,k} \diam(\G)$, 
  and for $k>1$,
  \begin{equation}
    \label{eq:2}
    \frac{\lambda_k}{\lambda_1} \leq 
     \frac{B_{k,k}}{B_{1,1}} \, \diam(\G) \; ,
  \end{equation}
where $\diam(\G)$ is the diameter of the network graph $\G$.
\end{proposition}

Proof of this result may be found in Appendix~\ref{sec:proof-bound-decay}.
The inequality in (\ref{eq:2}) indicates that the spectral decay 
in $G$ at worst parallels the decay of the edge betweenness on $\G$.
In fact, examination of these quantities for the Rocketfuel datasets
suggests that the decay in the $\lambda_k$ can be noticeably faster.
Nevertheless, Proposition~\ref{prop:1} provides a connection through
which recent and ongoing work on betweenness in complex networks
(\eg, \cite{barthelemy04:between,newman04:community}) 
may be found to have direct implications on the present context.

\section{Prediction of End-to-End Network Properties}
\label{sec:prediction-E2E-properties}

In this paper, we take as our monitoring goal the task of
obtaining accurate (approximate) knowledge of a linear summary of
network path conditions. That is, we seek to accurately predict a
linear function of the path conditions $y$, of the form $l^Ty$
where $l\in \R^{n_p}$, based on measurements, say $y_s\in\R^k$,
of a subset of $k$ paths.  Two
such linear summaries are the network-wide average, given by $l^T
y$ for $l_i \equiv 1/n_p$, and the difference between the averages
over two groups of paths $\P_1$ and $\P_2$, given by $l_i =
1/\abs{\P_1}$ for $i \in \P_1$ and $l_i = -1/\abs{\P_2}$ for $i \in
\P_2$.  The prediction of $l^Ty$ from the $k$ sampled path values
in $y_s$ can be viewed as a particular instance of the classical
problem of prediction in the statistical literature on sampling
\cite{valian00:finite_pop_sampling}. In this section, we (i) lay
out our statistical framework, (ii) describe an accompanying path selection
algorithm, and (iii) provide an analytical characterization 
of expected performance properties for our overall prediction methodology.

\subsection{Statistical prediction from sampled paths.}
\label{sec:stat-pred}

We begin by building a model for the end-to-end properties in $y$.
In the work that follows, it is only necessary that the first two
moments of $x$ and $y$ be specified, as opposed to a full
distributional specification.  Let $\mu$ be the mean of $x$ and
let $\Sigma$ be the covariance of $x$.  Then the corresponding
statistics for $y$ are simply $\nu=G \mu$ and $V = G \Sigma G^T$,
respectively.

Now fix $k \leq \rank(G)$.  Let $y_s\in\R^k$ denote the values
$y_{i_1},\dots, y_{i_k}$ of the metric of interest for $k$ paths
$i_1,\dots,i_k\in\P$ that are to be sampled (\ie,
measured), and let $y_r\in\R^{n_p-k}$ denote the values for those
$n_p-k$ paths that remain.  Similarly, let $G_s$ be those rows of
$G$ corresponding to the $k$ paths, $i_1,\dots,i_k$ and let $G_r$
be the remaining rows.  We have thus partitioned $y$ and $G$ into
$y = [y_s^T, y_r^T]^T$ and $G = [\Gs^T, \Gr^T]^T$, and we may
similarly re-express the mean and covariance of $y$ as
\begin{equation}
  \label{eq:3}
  \nu = 
  \begin{bmatrix}
    \nu_s \\ \nu_r
  \end{bmatrix}
  = 
  \begin{bmatrix}
    G_s \mu \\ G_r \mu
  \end{bmatrix}
\quad\text{and}\quad
  V = 
  \begin{bmatrix}
    \Vss & \Vsr \\
    \Vrs & \Vrr
  \end{bmatrix}
  = 
  \begin{bmatrix}
    \Gs \Sigma \Gs^T & \Gs \Sigma \Gr^T \\
    \Gr \Sigma \Gs^T & \Gr \Sigma \Gr^T
  \end{bmatrix}
  \text{.}
\end{equation}

If the standard mean-squared prediction error (MSPE) is used to
judge the quality of a predictor, \ie, if the quality of a
predictor $p(y_s)$ is measured by $\mspe(p(y_s))\equiv E[(l^Ty - p(y_s))^2]$, 
then the best predictor is known to be given by the conditional expectation
$E[l^Ty|y_s] = l_s^T y_s + E[l_r^Ty_r|y_s]$,
where 
$l=[l_s^T, l_r^T]^T$ is partitioned in the same manner as $y$.
But this predictor requires knowledge of the joint distributional
structure. It is therefore common practice to restrict attention
to a smaller and simpler subclass of predictors.  A natural choice
is the class of linear predictors, in which case, the best linear
predictor (BLP) is given by the expression
\begin{equation}
  \label{eq:4}
  a^T y_s = l_s^T y_s + l_r^T \Gr \mu + l_r^T c_* (y_s - \Gs \mu) \text{,}
\end{equation}
where $c_*$ is any solution to $c_* \Vss = \Vrs$.
However, without knowledge of $\mu$, the BLP in (\ref{eq:4}) is an
ideal that cannot be computed.  One natural solution is to
estimate $\mu$ from the data.
Using generalized least squares, the mean can be estimated as
$\muhat = [\Gs^T \Vss^{-1} \Gs]^{-} \Gs^T
\Vss^{-1} y_s$, where $M^{-}$ denotes a generalized inverse of a
matrix $M$.  Substituting $\muhat$ for $\mu$ in (\ref{eq:4})
produces an estimate of the BLP (an E-BLP) that is a function of
only the measurements $y_s$, the routing matrix $G$ and the link
covariance matrix $\Sigma$.  Specifically, we obtain the predictor
\begin{equation}
  \label{eq:5}
    \ahat^T y_s = l_s^T y_s + l_r^T \Gr [\Gs^T \Vss^{-1} \Gs]^{-} \Gs^T
    \Vss^{-1} y_s 
    = l_s^T y_s  + l_r^T \Vrs \Vss^{-1} y_s \text{.}
\end{equation}
The derivation of these and the other expressions above parallels
that of similar linear prediction methods in spatial
statistics---so-called `kriging' methods---which motivates
the name `network kriging'.  An example of the basic underlying argument,
in the case of simple linear statistical models, can be found in
\cite[pp.\ 225--227]{christensen87:plane}.  In the case of the present
context, the derivation requires only that $\Vss$ be invertible and 
that $\Sigma$ be positive definite.  

\subsection{Path Selection Algorithm}
\label{sec:path-selection}

The material in Section~\ref{sec:stat-pred} assumes a set of
measurements from $k$ paths $i_1,\ldots,i_k\in\mathcal{P}$.
However, given the resources to measure any $k$ paths in a
network, we are still faced with the question of which $k$ paths
to measure.  A natural response would be to choose $k$ paths that
minimize $\mspe(\ahat^T y_s)$, over all subsets of $k$ paths.
Standard manipulations yield that this quantity has the form
\begin{equation}
  \label{eq:6}
      \mspe(\ahat^T y_s) = 
      \underbrace{l_r^T \left( \Vrr -  \Vrs \Vss^{-1} \Vsr\right)
        l_r}_{\mspe(a^T y_s)} 
 + \underbrace{l_r^T \left(\Vrs \Vss^{-1} G_s - \Gr\right)
   \mu}_{(\bias \ahat^T y_s)^2}. 
\end{equation}
Of course, since we typically do not know $\mu$, minimization of
the full expression for $\mspe(\ahat^T y_s)$
is an unrealistic goal in practice.  Instead, if adequate
information on the covariance matrix $\Sigma$ is available, one
might consider trying to minimize $\mspe(a^T y_s)$.  A useful
equivalent expression for this quantity is
\begin{equation}
  \label{eq:7}
    \mspe(a^T y_s) 
    = l_r^T (\Gr C)
    (I - B_s)
    (\Gr C)^T l_r \,
        \text{,}
\end{equation}
where $C$ is a nonsingular matrix satisfying $\Sigma = C C^T$,
such as $\Sigma^{-\frac{1}{2}}$, and 
$B_s$ is the orthogonal projection matrix onto the span of the rows
of $\Gs C$, \ie, onto $\row(\Gs C)$.  Since orthogonal projection matrices
are idempotent and symmetric, the MSPE in (\ref{eq:7}) can be viewed as 
the square of the Euclidean norm of the projection of $(\Gr C)^T l_r$
onto the complement of $\row(\Gs C)$, \ie, onto 
$\row(\Gs C)^\perp = \nul(\Gs C)$.

In order to better appreciate the interpretation of (\ref{eq:7}), consider
the special case of predicting a single unmeasured path (\ie, 
$l\equiv 0$ except for a single $1$ in $l_r$), with $\Sigma=I$.
The MSPE in (\ref{eq:7}) then simply measures the extent to which the
corresponding row of $G$ for this path lies outside of $\row(\Gs)$.
Similarly, the more interesting case of a non-trivial $l$ can be
interpreted roughly as seeking a subset of $k$ paths for whom the
rows in $G_s$ capture as many of the rows of $G$ as possible to
the largest extent possible.

From the standpoint of optimization theory, our path-selection
problem may be viewed as an example of the so-called `subset
selection' problem in computational linear algebra.  In the case
just described, and more generally for diagonal $\Sigma$, the
selection of an appropriate subset of rows of $GC$ has a
meaningful physical interpretation, in terms of the selection of
paths, and vice versa.  Exact solutions to this problem are
computationally infeasible (it is known to be NP-complete), but
the problem is well-studied and an assortment of methods for
calculating approximate solutions abound.

The method we have used for the empirical work in this paper was
adapted from the subset selection method described in Algorithm
12.1.1 of \cite{golub89}.  Essentially, our algorithm
makes heuristic use of a QR-factorization with column pivoting to
find $k$ rows of $G$ that approximate the span of the first $k$
left singular vectors of $GC$. The left singular vectors form an
orthonormal basis for the range of $GC$ and the magnitude of their
corresponding singular values indicates their relative importance.
Note that these singular values are precisely the square-root of
the eigenvalues of $(GC)^T(GC)$, \ie, the decaying spectrum from
Section~\ref{sec:routing-matrices}. See \cite{chua:framework} for
additional details.

For a given choice of $k$, the overall complexity for the
computation of the E-BLP in (\ref{eq:5}) is dominated by the
computation of the SVD of $GC$, which is $O(n_p^2 n_e)$.  This can
likely be improved through the use of methods for sparse matrices,
since the entries of $GC$ tend to include a large fraction of
zeros.  The other components of the computation 
are the QR-factorization with column
pivoting, which is $O(k^2 n_p)$, and the computation of $\Vss^{-1}$
which is only $O(k^3)$.

\subsection{Characterization of MSPE Properties}
\label{sec:properties}

Analytical arguments are useful for characterizing the expected 
performance of a predictor resulting from the combination of 
equation (\ref{eq:5}) with an arbitrary subset selection algorithm.
For example, we have the following bound.
\begin{proposition}
  \label{prop:2}
  Denoting the $i^{\text{th}}$ row of $G$ as $G_{(i)}$\,, let
  $p_i = \norm{G_{(i)}}_2^2 / \norm{G}_F^2$\,, where
  $\norm{\cdot}_2$ and $\norm{\cdot}_F$ are the matrix 2-norm and
  the Frobenius matrix norm, respectively.  Let $\tildeGs$ be a
  rescaled version of $\Gs$\,, under the operation $G_{(i)}
  \rightarrow G_{(i)}/\sqrt{k p_i}$\,, for each of the $k$ rows in
  $G_s$.  Then if
  \begin{equation}
    \label{eq:8}
    \Norm{G^TG - \tildeGs^T\tildeGs}_F \leq f(k) \Norm{G}_F^2\enskip ,
  \end{equation}
for some $f(k)$, the MSPE can be bounded as
  \begin{equation}
    \label{eq:9}
    \mspe(\ahat^T y_s) 
    \leq (\Norm{\mu}_2^2 + 1)
    \left(\lambda_{k+1} + 2 f(k) \Norm{G}_F^2 \right)
    \Norm{l}_2^2 \text{.}
  \end{equation}
\end{proposition}

Proof of this result may be found in Appendix~\ref{sec:proof-mspe-char}.
The inequality in (\ref{eq:9}) shows that the decay of the MSPE 
in $k$ is controlled by two factors.  The first factor, $\lambda_{k+1}$,
quantifies how well $G$ may be approximated by
its first $k$ singular dimensions, and will be small when $k$ is no
smaller than the effective rank of $G$.  The second factor, $f(k)$,
quantifies the ability of the underlying subset selection algorithm
to approximate $G$ by a matrix $G_s$ formed from $k$ of its rows.
Our empirical experience indicates that 
in practice $\lambda_{k+1}$ can be much smaller than the term
involving $f(k)$, which suggests that the rate of decay
is dominated by $f(k)$.

\begin{figure}[tbp]
  \centering
  \includegraphics[width=\onefig]{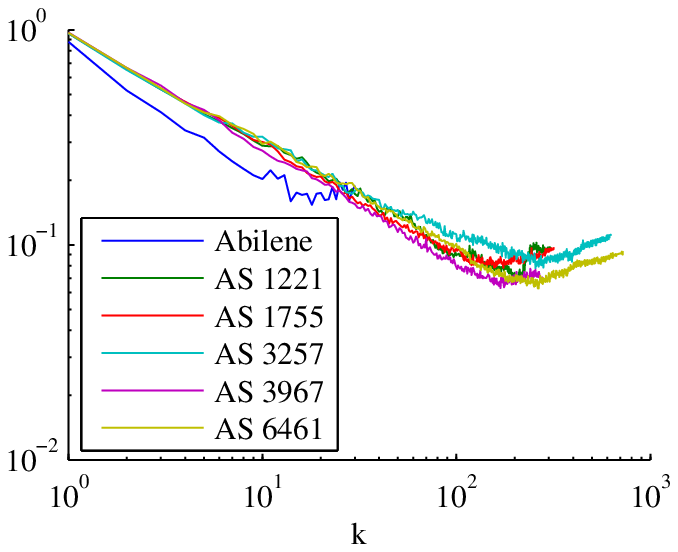}
  \caption{Plot of \textsf{$\Norm{G^TG - \tildeGs^T\tildeGs}_F / \norm{G}_F^2$} for the Abilene
    backbone and Rocketfuel
    networks.} 
  \label{fig:fk}
\end{figure}

To get an idea of the behavior of $f(k)$ for real networks, we
computed $\Norm{G^TG - \tildeGs^T\tildeGs}_F / \norm{G}_F^2$ for the
Abilene backbone and five networks mapped by the Rocketfuel project, for the
subset selection algorithm described in Section~\ref{sec:path-selection}.
As can be seen in Figure~\ref{fig:fk}, there is a strong power-law decay in
$\Norm{G^TG - \tildeGs^T\tildeGs}_F / \norm{G}_F^2$, with an
exponent that ranges from $-0.49$ to $-0.53$.  The consistency in
this decay is quite remarkable, as it says that our subset selection
algorithm finds a similarly good set of paths, for each $k$, in each
of the networks.
That this decay in $f(k)$ indeed translates into good MSPE properties
will be seen in Section~\ref{sec:empirical-evaluation}.

Proposition~\ref{prop:2} holds for any given deterministic path
selection algorithm.  Perhaps surprisingly, it is possible to 
state a similar result for a randomized path selection algorithm.
The empirical work we present in the
next section clearly establishes the practical effectiveness of
our methodology using the deterministic algorithm described in
Section~\ref{sec:path-selection}.  However, effective randomized
path selection algorithms are interesting to consider in the sense that
such algorithms would be able to reduce the likelihood that small, highly
localized events consistently remain outside the span of the
sampled paths.  The following result, the proof of which follows
as a direct corollary of Theorem 3 in~\cite{drineas04}, suggests
that an algorithm that randomly selects paths roughly in
proportion to path length could indeed achieve similar performance
characteristics.
\begin{proposition}
  \label{prop:3}
  Let $\tildeGs$ be a matrix constructed as in
  Proposition~\ref{prop:2}, but now from $c$ paths randomly selected
  (with replacement) with respect to the probabilities 
  $\{p_i\}$.  If the estimator $\ahat^T y_s$ in (\ref{eq:5}) is
  constructed based on the first (at most) $k\le c$ singular dimensions
  of $\tildeGs$, then for any $c\leq n_p$ and $\delta > 0$,
  \begin{equation}
    \label{eq:10}
    MSPE(\ahat^T y_s) \leq (\norm{\mu}_2^2 + 1) 
    (\lambda_{k+1} + 2(1+\sqrt{\ln(2/\delta)}) c^{-\frac{1}{2}} \,
    \norm{G}_F^2) \, \norm{l}_2^2
  \end{equation}
  holds with probability at least $1-\delta$. 
\end{proposition}

\section{Empirical Validation of Prediction Methodology}
\label{sec:empirical-evaluation}

In this section, we show how our framework may be applied to
address two practical problems of interest to network providers
and customers.  In particular, we show how the appropriate
selection of small sets of path measurements can be used to
(i) accurately estimate network-wide averages of path delays, and (ii)
reliably detect network delay anomalies.  We first describe the 
assembly of our dataset, and then present the two applications.

\subsection{Data: Abilene Path Delays}
\label{sec:data-abilene}

Our methods are applicable to any per-link metric that adds to
form per-path metrics.  As a specific example, we consider delay.
In order to validate our prediction methodology, we constructed a full set
of path-delay data for the Abilene network, using measurements
obtained from the NLANR Active Measurement Project (AMP).  
This project continually performs traceroutes between
all pairs of AMP monitors on ten minute intervals.  Because most
AMP monitors are on networks with Abilene connections, most
traceroutes pass over Abilene.  These data thus provide a highly detailed
view of the state of the Abilene network.

Beginning with a full set of measurements taken over 3 days in 2003,
our approach to constructing per-path delays from this data
consisted of (i) estimating per-link delays $x^{(t)}\in\R^{30}$
across the 30 network links, for each consecutive ten minute epoch $t$,
and (ii) computing the per-path delays from these per-link delays, on 
the 110 Abilene paths, using the relation
$y^{(t)}=Gx^{(t)}$, where $G$ is a fixed routing matrix corresponding 
to the routing on Abilene at the start of the 3 day period.
The end result is a temporally indexed sequence of path-delay vectors
$y^{(t)}$ over $432$ successive epochs during the 3 day period.
Note that the inferred link delays $x^{(t)}$ are not explicitly used by our 
prediction methodology, but rather were only necessary for the construction 
of the path delays, after which they were discarded.

To construct the link delays, for each epoch $t$,
we started with traceroutes between the 14,917
pairs of AMP monitors for which complete data were available.  
Links comprising the
Abilene network were identified by their known interface
addresses.  Since different traceroutes traverse each link at
slightly different times, and since each traceroute takes up to
three measurements per hop, we formed a single estimate of each
link's delay by averaging across all the traceroutes that measured
that link in the current epoch.  This approach yields a single measure of
delay for each link and epoch; while this measure does not capture
the variations in delay that occur within a ten minute interval,
it provides a realistic and representative value for delay that
is sufficient for our validation purposes.

Recall that the statistical framework in Section~\ref{sec:stat-pred}
involved only the first two moments of the link delays, the
mean $\mu$ and the covariance $\Sigma$.  While our methodology
does not require knowledge of $\mu$, since it is estimated at each epoch 
interval as part of the calculation of our predictor, we note here that the 
mean delays in our data for Abilene's $30$ directed edges were fairly 
uniform over the interval from 2 to 36 milliseconds, with standard
deviations that ran from 0.16 to 0.94 over the full three-day
period.  Our methodology does, however, require knowledge of $\Sigma$,
which in practice must be elicited from either historical data
or possibly periodic, infrequent measurements on the links.  For the
purposes of the validation in this section, we used the per-link 
delays $x^{(t)}$ for one day's worth of data to obtain an estimate of $\Sigma$.
The entries in this matrix were found to be primarily dominated by the 
diagonal elements, with only a small number of off-diagonal entries
of similar magnitude.  Inspection of the actual delay data suggested
that the majority of the latter were due to artifacts of the measurement 
procedure.  Therefore, in implementing our methodology we
took $\Sigma$ to be the corresponding diagonal matrix.
See~\cite{chua:framework} for additional details, including
evidence regarding the improvement gained over the choice $\Sigma\propto I$.

\subsection{Monitoring a Network-Wide Average}
\label{sec:network-wide-average}

\begin{figure}[tbp]
  \centering
  \includegraphics[width=\onefig]{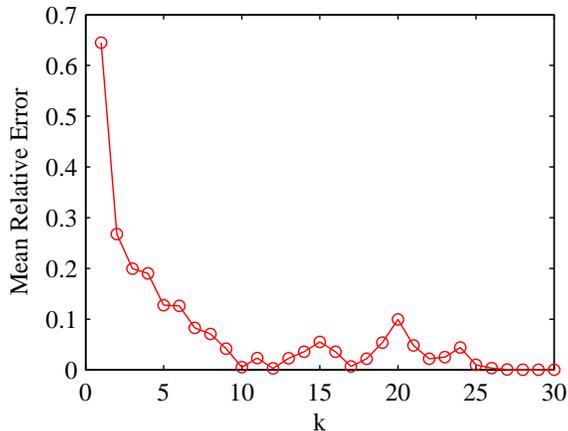}
  \caption{Mean relative prediction error as a function of $k$.} 
  \label{fig:rel-error}
\end{figure}

An average is perhaps the most basic network-wide quantity that
one might be interested in monitoring. So as our first
application, we consider the prediction of the average delay over
all $n_p=110$ Abilene paths as a function of time $t$, \ie,
$l^Ty^{(t)}$ with $l_i\equiv 1/n_p$, $i=1,\ldots,n_p$, and 
$t=1,\ldots,432$.
Using (\ref{eq:5}), we computed
predictions of the network-wide average path delay during each
epoch, for a choice of $k=1,\ldots,30$ measured paths. The paths
were chosen using the algorithm described in
Section~\ref{sec:path-selection}. 
To summarize the accuracy of our
predictions, we calculated the average relative error for each
$k$, where the average is taken all the 432 epochs.  The results are
shown in Figure~\ref{fig:rel-error}.

Recall that exact recovery of all paths delays $y^{(t)}$, at a
given epoch $t$, requires measurement of $k^*=\rank(G)$ paths,
which in this case means $k^*=30$.  In examining Figure~\ref{fig:rel-error},
note that in comparison a relative error of roughly only $10\%$ is
achievable using only $k=7$ path measurements.
Increasing $k$ further improves
the accuracy of the prediction up until around $k=9$ or $10$,
after which it basically levels out.  Since it is roughly at this
point that the spectra of the (weighted) routing matrices level
out as well, this suggests that our subset selection algorithm is
indeed doing what we are asking of it, in that it is tracking the
effective rank quite closely.  Additional results of this sort, on
a variety of simulated datasets, may be found 
in~\cite{chua05:efficient}.

\begin{figure}[tbp]
  \centering
  \includegraphics[width=\onefig]{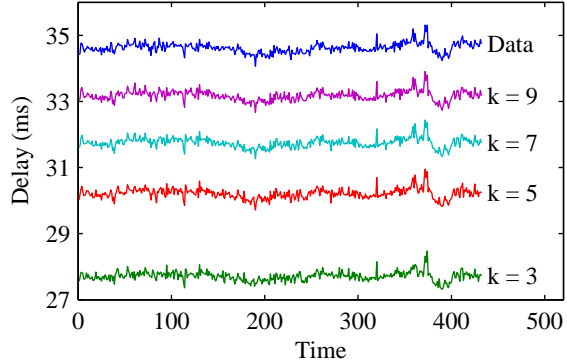}
  \caption{Predictions of network-wide average path delays, for
           various choices of $k$.}
  \label{fig:GC-preds}
\end{figure}

To get a better idea of how well the predictors performed, we can
compare plots of the predictions against a plot of the actual mean
delays, as shown in Figure~\ref{fig:GC-preds} for $k=3,5,7$ and
$9$.  Note that all of the predictions mirror the rise and fall of
the actual network-wide delay quite closely---even for $k=3$
measured paths the correlation with the data time-series is
$\rho=0.814$. However, it is also clear that there is a downward
bias in these predictions, and that this bias is increasingly
prominent as $k$ decreases.
The source of this bias can be traced to a lack of information on
links in the network that are traversed by none of the $k$
measured paths.  In fact, the generalized inverse used in our
E-BLP in (\ref{eq:5}) simply estimates the corresponding values
$x_j$ on these links to be zero.  Hence, as we reach a point where
every link contributes to at least one measured path, as it does
by roughly $k=10$, the bias diminishes accordingly.  Note,
however, that the bias for each $k$ in Figure~\ref{fig:GC-preds}
is fairly constant.  This suggests that a small amount of
additional measurement information could go a long way.

\begin{figure}[tbp]
  \centering
  \includegraphics[width=\onefig]{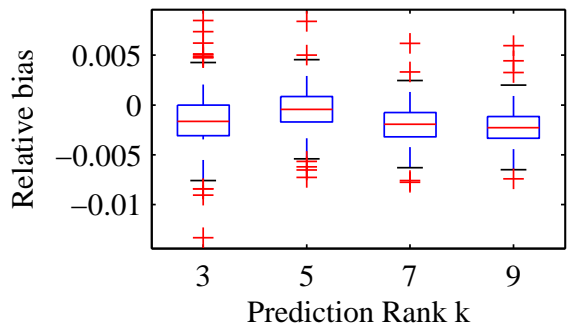}
  \caption{Relative bias after bias correction.}
  \label{fig:bias-corrected-ts}
\end{figure}

We implemented a simple method of bias correction, that uses a
one-time measurement of a sufficient set of paths for complete
reconstruction of the link delays (in this case, 30 paths). Since
it is a one-time-only measurement, it represents a minimal
addition to the network measurement load.  The bias of our
prediction for the first epoch was then calculated and used to
adjust the predictions in the other 431 epochs, which amounts to
a simple shift upward of each curve in Figure~\ref{fig:GC-preds}.
Boxplots of the relative bias remaining after application of this
procedure are shown in Figure~\ref{fig:bias-corrected-ts}.  The
predictions are now extremely accurate, usually being off by less
than $0.3\%$, and almost always within $1\%$---even when as few as
$k=3$ paths are measured.

Before moving on, we note that this performance on whole networks
extends to subnetworks. In particular, we have successfully used our
methods to make reliable comparisons between sub-network delays
in the context of multi-homing~\cite{chua:framework}.

\subsection{Anomaly Detection}
\label{sec:anomaly-detection}

The application in Section~\ref{sec:network-wide-average}
evaluates our predictor by standard statistical summaries, in
essence looking at the accuracy of the predictor at hitting an
unknown target. But it is also important to evaluate the accuracy
in terms of accomplishing higher-level tasks. One such
higher-level task of importance is the detection of potentially
anomalous events.

\begin{figure}[tbp]
  \centering
  \includegraphics[width=\onefig]{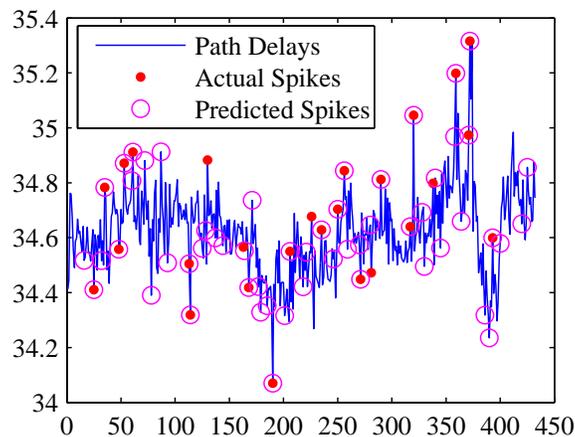}
  \caption{Comparison of predicted and actual spikes. The real spikes
  are those that exceed 3 times the standard deviation of the previous
  6 epochs. The predicted spikes are those epochs where the rank 9
  prediction exceeds 2 times the standard deviation of its previous 6
  epochs. }
  \label{fig:spikes}
\end{figure}

For the purposes of illustration with our Abilene delay data, 
we define an anomaly as a spike in the network-wide average
path delay that deviates from the mean of the
previous six values (\ie, one hour) by more than a
prescribed amount. 
For example, the dots in
Figure~\ref{fig:spikes} indicate points at which the average path
delay differs from the mean of the previous six epochs by more
than three times their standard deviation.

\begin{figure}[tbp]
  \centering
  \includegraphics[width=\onefig]{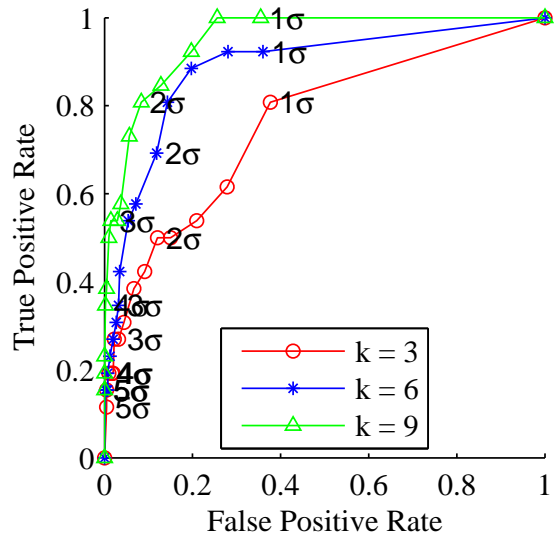}
  \caption{ROC curves for predicting $3\sigma$ spikes. The threshold
  used to predict the spikes is varied from $1\sigma$ 
    to $5\sigma$ in increments of $0.25\sigma$. }
  \label{fig:ROC-threshold}
\end{figure}

To predict when such anomalies occur, we look for spikes in the
predicted average path delay, calculated as described in
Section~\ref{sec:network-wide-average} and using a user-defined
threshold.  It is interesting to examine the effect of choice of
both $k$ and this threshold parameter.  Insight can be obtained
by examining ROC (Receiver Operating Characteristic)
curves such as those in Figure~\ref{fig:ROC-threshold}.  Such
plots, showing the true positive rate against the false positive rate
for different parameter values, are a common tool for establishing
cutoff values for detection tests.  Each curve in
Figure~\ref{fig:ROC-threshold} is formed by taking one value for
$k$ and varying the threshold level.  Examining these curves, one
sees that for a given threshold, say $1\sigma$, the true positive
rate increases with the sample size $k$ while the false positive
rate stays about the same. Working with a $k=9$ prediction, we see
that the upper-left corner of the ROC curve (the best trade-off
between a low false positive rate and a high true positive rate)
occurs at around $2\sigma$.

In Figure~\ref{fig:spikes}, the results are shown for the case
$k=9$, with a threshold of $2\sigma$.  Circles have been placed
along the actual path delay time series at the epochs that were
flagged as anomalies in the predicted time series.  On the whole,
this predictor is quite accurate. Most of the major spikes are
flagged, resulting in a true positive rate of 81\%, while the
false positive rate is only 8\%.  Furthermore, most of these false
positives seem to occur at lesser spikes in the actual delay data.

\section{Discussion}
\label{sec:discussion}

The identification of an inherent statistical prediction problem 
in the task of end-to-end network path monitoring---which we
have dubbed `network kriging'---is analogous in its potential impact with
the identification of, say, traffic matrix estimation with tomography
(\ie, `network tomography').  In both cases, the identification 
serves as a critical pointer to an already established literature, through
which methodology may be developed by leveraging various principles
and tools.  Here we have focused exclusively on one basic version of
the network kriging problem, the prediction of linear metrics
of path properties using linear modeling principles.  However, there is
ample room for work beyond this, including for example extensions to
non-linear metrics and temporal prediction models.

We have successfully demonstrated the promise of our proposed
methodology on data obtained from a real network.  Nevertheless,
there are various practical issues to be dealt with to optimize
our framework for full-scale implementation.  These include
efficient strategies for dealing with monitor failures, link
failures and routing changes.  We expect that
many of these may be addressed using tactics similar to those
proposed in~\cite{chen04:algebraic}.

On a final note, we mention again the fundamental importance of the
material in Section~\ref{sec:reduced-rank}, on the prevalence of
low effective rank
of routing matrices, to the success of our methodology.  Put simply,
low effective rank allows for the possibility of
effective network-wide monitoring with reduced measurement sets.  While
we have exploited this characteristic for the purpose of path monitoring,
it should apply equally well when the goal involves selective monitoring
of links. The driving factors responsible for the low effective
rank are poorly understood and remain an interesting open problem. 

\appendices

\section{Proof of Proposition~\ref{prop:1}.}
\label{sec:proof-bound-decay}
Let $\lambda_1 \geq \dots \geq \lambda_n$ denote the
eigenvalues of the matrix $B$
and define the spectral
radius of $B=G^TG$ as
$  \rho(B) = \max\{ \abs{\lambda_1},\dots, \abs{\lambda_n} \} = \lambda_1$.
Now partition the matrix $B$ into 
$B =
\bigl[
\begin{smallmatrix}
  E & C^T \\
  C & F^{\phantom{T}}
\end{smallmatrix} 
\bigr]$ where $F$ is an $m\times m$ matrix with $m < n_e$ whose
diagonal elements are the $m$ smallest betweennesses.  Denote the
eigenvalues of $F$ by $\theta_1 \geq \dots \geq \theta_{m}$ and,
for convenience, let $\theta_i=-\infty$ for $i>m$ and
$\theta_i=\infty$ for $i\leq0$.  Cauchy's Interlace Theorem tells
us that the $m$ eigenvalues of $F$ are upper bounds for the $m$
smallest eigenvalues of $B$, respectively.

For the moment, we focus on the case of $\theta_1
\geq\lambda_{n-m+1}$. If we can bound $\theta_1$ we will have a
bound for $\lambda_{n-m+1}$.  Define the $i$-th row sum of
$B$ as $r_i(B) = \sum_{j=1}^n \abs{B_{i,j}}$, and the deleted row
sum of $B$ as $\tilde{r}_i(B) = \sum_{\substack{j=1 \\ j \neq
    i}}^n \abs{B_{i,j}}$.  It then follows, by way of a corollary
to Gershgorin's Theorem\cite{varga04:gersgorin}, that the spectral
radius of $F$ is bounded by
\begin{equation}
  \label{eq:12}
  \rho(F) \leq \max_{1 \leq i \leq m} r_i(F) 
  \leq \max_{n-m+1 \leq i \leq n} r_i(B) \; .
\end{equation}
Each of these row sums $r_i(B)$ can be bounded in terms of the
diagonal element $B_{i,i}$. To see this, note that each of the
$B_{i,i}$ paths that use edge $i$ have a length of at most
$\diam(\G)$.  So each path can contribute at most $\diam(\G)$ to the
row sum $r_i(B)$.  This means that for a betweenness matrix $B$ we
have $r_i(B) \leq B_{i,i} \, \diam(\G)$.  Recalling that we have
ordered the edges such that the diagonal elements $B_{i,i}$ are
non-increasing, our bound for the spectral radius becomes
\begin{equation}
  \label{eq:13}
  \rho(F) \leq B_{n-m+1,n-m+1} \diam(\G) \; .
\end{equation}
Which means that our bound for the eigenvalue of $B$ is
\begin{equation}
  \label{eq:14}
  \lambda_{n-m+1} \leq \theta_1 \leq r_{n-m+1}(B) \leq B_{n-m+1,n-m+1} \, \diam(\G)
  \, .
\end{equation}
Finally, note that we are free to choose
$m=1,\dots,n-1$. So for $i = 2, \dots, n$ we have an upper bound that
decays no worse than the betweenness $B_{i,i}$, namely,
\begin{equation}
  \label{eq:15}
  \lambda_i \leq B_{i,i} \diam(\G)\;.
\end{equation}
The $i=1$ case follows from Gershgorin applied directly to $B$.

To establish the second bound, given in \eqref{eq:2} of
Proposition~\ref{prop:1}, note that for any unit vector $u$ we
have $\norm{Bu}_2 \leq \lambda_1$.  In particular, for $u = e_1 =
[1,0,\dots,0]^T$ we have $\norm{B e_1}_2 \leq \lambda_1$.  Since
$B e_1$ is the first column of $B$, we have $\norm{Be_1}_\infty =
B_{1,1}$. Thus $B_{1,1} = \norm{B e_1}_\infty \leq \norm{B e_1}_2
= \lambda_1$.  Dividing the left and right-hand sides of
\eqref{eq:15} by $\lambda_1$ and $B_{1,1}$, respectively yields
the desired bound.

\section{Proof of Proposition~\ref{prop:2}.}
\label{sec:proof-mspe-char}
Since 
$\mspe(\ahat^T y_s) = \Norm{\mu^T (I-\Bs) G^T l}_2^2 
  + \Norm{(I-\Bs) G^T l}_2^2 
\leq (\Norm{\mu}_2^2 + 1)\, \norm{(I - \Bs)G^T}_2^2 
  \,\Norm{l}_2^2$, 
we need only show $\norm{(I - \Bs)G^T}_2^2 \leq \lambda_{k+1} + 2
f(k) \norm{G}_F^2$.  To do so, we proceed as in the proof of Theorem~3 in
\cite{drineas04}.  Note that for any vector $x\in \R^n$ with
$\Norm{x} = 1$ we can write $x = ay + bz$, where $y\in\Ss =
\row(\Gs)$, $z\in\Ss^\perp$, $a,b\in\R$, and $a^2+b^2 = 1$. Using
this decomposition of $x$ and the sub-linearity of the $2$-norm we
can write
\begin{align}
  \label{eq:20}
  \Norm{G^T - \Bs G^T}_2 &=\max_{\Norm{x}=1} \Norm{x^T(G^T - \Bs G^T)} \\
  &\leq \max_{\substack{y\in \Ss\\ \Norm{y}=1}} \Norm{y^T (G^T - \Bs G^T)} 
  + \max_{\substack{z\in \Ss^\perp\\ \Norm{z}=1}} \Norm{z^T (G^T - \Bs G^T)}
\end{align}
At this point, we note that since $y\in\Ss$ we have $y^T \Bs=y^T$,
which means that $y^T(G^T - \Bs G^T) = 0$. Furthermore,
$z\in\Ss^{\perp}$ means that $z^T \Bs=0$.  So the
upper bound becomes
\begin{equation}
  \label{eq:21}
  \Norm{G^T - \Bs G^T}_2 \leq
  \max_{\substack{z\in\Ss^{\perp}\\\Norm{z}=1}} \Norm{z^T G^T}.
\end{equation}
We can bound the maximum in \eqref{eq:21} in terms of $\tildeGs$ by
\begin{equation}
  \label{eq:22}
  \norm{z^T G^T}_2^2 
  = z^T(G^T G - \tildeGs^T\tildeGs) z + z^T \tildeGs^T\tildeGs z 
  \leq \norm{G^TG - \tildeGs^T\tildeGs}_F + \sigma_{k+1}^2(\tildeGs) \text{,}
\end{equation}
where $\sigma^2_{k+1}(M) \equiv \lambda_{k+1}(M)$,
since $\norm{z}_2 = 1$ and $\norm{M}_2 \leq
\norm{M}_F$. Combining \eqref{eq:21} and \eqref{eq:22} gives us
\begin{equation}
  \label{eq:23}
  \Norm{G^T - \Bs G^T}_2^2 \leq \sigma_{k+1}^2(\tildeGs) + \Norm{G^TG -\tildeGs^T\tildeGs}_F.
\end{equation}
Turning to Corollary~8.6.2 in \cite[p.~449]{golub89} we have for
$k=1,\dots,n_e$
\begin{equation}
  \label{eq:24}
  \abs{\sigma_{k+1}(G^TG) - \sigma_{k+1}(\tildeGs^T\tildeGs)} 
  \leq \norm{G^TG - \tildeGs^T\tildeGs}_2 .
\end{equation}
Using the bound for $\norm{G^TG - \tildeGs^T\tildeGs}_F$ that we
are given in \eqref{eq:8}, and the relation $\norm{M}_2 \leq
\norm{M}_F$ for any matrix $M$ we have
\begin{equation}
  \label{eq:25}
  \norm{G^TG - \tildeGs^T\tildeGs}_2 \leq f(k)\, \norm{G}_F^2 .
\end{equation}
Thus $
  \abs{\sigma_{k+1}(G^TG) - \sigma_{k+1}(\tildeGs^T\tildeGs)} 
  = \abs{\sigma_{k+1}^2(G) - \sigma_{k+1}^2(\tildeGs)} 
  \leq \Norm{G^TG - \tildeGs^T\tildeGs}_2 
  \leq f(k)\, \norm{G}_F^2$,
which leads us to
$  \sigma_{k+1}^2(\tildeGs) \leq f(k) \Norm{G}_F^2 + \sigma_{k+1}^2(G)$.
Combined with \eqref{eq:23} and \eqref{eq:25} this yields that
$\Norm{(I-\Bs)G^T}_2^2\leq \lambda_{k+1} + 2 f(k) \Norm{G}_F^2$,
as was to be shown.

\bibliographystyle{IEEEtran}
\bibliography{JSAC}


\end{document}